\documentclass[11pt]{amsart}

\usepackage{amssymb,amsfonts,latexsym}
\usepackage[centertags]{amsmath}
\usepackage{amsfonts}
\usepackage{amssymb}
\usepackage{amsthm}
\usepackage[all]{xy}
\usepackage{cite}

\usepackage{todo}

\newtheorem*{Prob*}{Problem}
\newtheorem{cont}{cont}[section]
\newtheorem{teo}[cont]{Theorem}
\newtheorem*{teo*}{Theorem}
\newtheorem{prop}[cont]{Proposition}
\newtheorem{lemma}[cont]{Lemma}
\newtheorem{cor}[cont]{Corollary}
\newtheorem{defn}[cont]{Definition}
\newtheorem*{Enunciato*}{Enunciato}

\numberwithin{equation}{section}

\theoremstyle{remark}
\newtheorem{remark}[cont]{Remark}

\newtheorem*{not*}{Notation}

\newcommand{\shF}{\mathcal{F}}

\newcommand{\shG}{\mathcal{G}}

\newcommand{\shH}{\mathcal{H}}
\newcommand{\shI}{\mathcal{I}}

\newcommand{\shE}{\mathcal{E}}
\newcommand{\PP}{\mathbb{P}}

\newcommand{\OO}{\mathcal{O}}
\newcommand{\odi}[1]{\mathcal{O}_{#1}}
\newcommand{\arr}{\longrightarrow}

\DeclareMathOperator{\Hl}{H} \DeclareMathOperator{\h}{h}
 
\DeclareMathOperator{\rk}{rk} \DeclareMathOperator{\Hom}{Hom}
 \DeclareMathOperator{\Aut}{Aut}
\DeclareMathOperator{\de}{deg}

\DeclareMathOperator{\depth}{depth}
 \DeclareMathOperator{\di}{dim}

\DeclareMathOperator{\Ext}{Ext}

\begin{document}
\title{n-dimensional Fano varieties of wild representation type}

\author[Rosa M. Mir\'{o}-Roig]{Rosa M. Mir\'{o}-Roig$^{*}$}
\address{Facultat de Matem\`atiques, Department d'Algebra i Geometria, Gran Via des les Corts Catalanes 585, 08007 Barcelona, Spain}
\email{miro@ub.edu}

\author[Joan Pons-Llopis]{Joan Pons-Llopis$^{**}$}
\address{Facultat de Matem\`atiques, Department d'Algebra i Geometria, Gran Via des les Corts Catalanes 585, 08007 Barcelona, Spain}
\address{Facultad de Matem\'{a}ticas, Departamento de \'{A}lgebra, Plaza de las Ciencias 3, 28040 Madrid, Spain}
\email{jfpons@ub.edu}

\date{\today}
\thanks{$^*$ Partially supported by MTM2010-15256.\\
$^{**}$ Supported by the research project MTM2009-06964.}

\begin{abstract} The aim of this work  is to provide the first examples of $n$-dimensional varieties of wild representation type, for arbitrary  $n\geq 2$. More precisely, we prove that all Fano blow-ups of $\PP^n$ at a finite number of points are of wild representation type exhibiting families
of dimension of order $r^2$ of simple (hence, indecomposable) ACM rank $r$ vector bundles for any $r\geq n$. In the two dimensional case, the vector bundles that we construct are moreover Ulrich bundles and $\mu$-stable with respect to certain ample divisor.
\end{abstract}

\maketitle

\tableofcontents

\section{Introduction}
A possible way to measure the complexity of a given ACM $n$-dimensional projective variety $(X, \OO_X(1))$ is to ask for the families of non-isomorphic indecomposable ACM vector bundles that it supports (recall that a vector bundle $\shE$ on $X$ is ACM if $\Hl^i(X,\shE\otimes \OO_X(t))=0$ for all $t\in\mathbb{Z}$ and $i=1,\dots,n-1$). The first result on this direction was Horrocks' theorem which states that on $\PP_k^n$ the only indecomposable ACM bundle up to twist is the structure sheaf $\OO_{\PP_k^n}$. Later on,  Kn\"{o}rrer in \cite{Kno} proved that on a smooth quadric hypersurface X the only indecomposable ACM bundles up to twist are $\odi{X}$ and the spinor bundles $S$.

Varieties that admit only a finite number of indecomposable ACM bundles (up to twist and isomorphism) are called of \emph{finite representation type} (see \cite{DG} and references herein). These varieties are completely classified (cfr.\cite{BGS} and \cite{EH}): They are either a projective space $\PP_k^n$, a smooth quadric hypersurface $X\subset \PP^n$, a cubic scroll in $\PP_k^4$, the Veronese surface in $\PP_k^5$ or a rational normal curve.

\vskip 2mm
If we look at the other extreme of complexity we would find the varieties of \emph{wild representation type}, namely, varieties for which there exist $r$-dimensional families of non-isomorphic rank $r$ indecomposable ACM bundles for arbitrary large $r$. In the case of dimension one, it is known that curves of wild representation type are exactly those of genus larger or equal than two.  Note that all vector bundles $\shE$ on a smooth curve $C$ are ACM  since the condition $\Hl^i(C,\shE(t))=0$ for all $t\in\mathbb{Z}$ and $1\le i \le \dim(C)-1$ is vacuous. For varieties of larger dimension, a complete classification of such varieties seems out of reach. As far as we know, in \cite{CH}, Casanellas and Hartshorne pointed out the first example of surface of wild representation type: the  smooth cubic surfaces in $\PP^3$. Unfortunately, the proof of the indecomposability of the Ulrich bundles that they construct  relies on  \cite[Theorem 5.3]{CH}, which proof is up to our knowledge incomplete. In \cite{PT}, it was shown that Del Pezzo surfaces of degree $\leq 6$ are of wild representation type surfaces.

\vskip 2mm
The main goal of this paper is to provide the first examples of $n$-dimensional varieties of wild representation type, for arbitrary  $n\geq 2$. Our source of examples will be Fano blow-ups $X=Bl_Z \PP^n$ of $\PP^n$ at a finite set of points $Z$. Recall that a smooth variety $X$ is called Fano if the anticanonical line bundle $-K_X$ is ample. In order the anticanonical line bundle $-K_X$ of a blow-up $X=Bl_Z \PP^n$ of $\PP^n$ at $s$ points to be ample, we are going to be force to impose some restrictions on $s$. Namely, in the $2$-dimensional case, i.e., for del Pezzo surfaces, it is classically known that we can take up to $8$ points. In higher dimensions, only the blow-up of up to one point is possible (cfr. \cite[Theorem 1]{Bon}). Our main result is the following one (see Theorem \ref{higheracm}):

\begin{teo*}
Let $X=Bl_Z \PP^n$ be a Fano blow-up of points in $\PP^n$, $n\geq 2$. Then, for any $r\geq n$, there exists a family of rank $r$ simple (hence, indecomposable)  ACM vector bundles of dimension $\frac{(n+2)n-4}{4}r^2-cr+c^2+1$ if $n$ is even or $((n+2)n-4)r^2-2cr+c^2+1$ if $n$ is odd. In particular, Fano blow-ups are of wild representation type.
\end{teo*}

In the $2$-dimensional case, i.e., for del Pezzo surfaces, much more information is obtained. In fact, the vector bundles that we construct share another particular feature: the associated module $\oplus_{t}\Hl^0(X,\shF(t))$ has the maximal number of generators. This property was isolated by Ulrich in \cite{Ulr}, an ever since modules with this property have been called Ulrich modules and correspondingly Ulrich bundles in the geometric case. The existence of rank $2$ Ulrich bundles on arbitrary del Pezzo surfaces were established by Eisenbud, Schreyer and Weyman in \cite{ESW}, Corollary 6.5.  We refer to \cite{ESW}, \cite{CH} and \cite{CKM} for more details on Ulrich bundles.
We show (see Theorem \ref{main}):

\begin{teo*}
Let $X$ be a del Pezzo surface of degree $d$ over an algebraically closed field $k$ of characteristic zero. Then, for any $r\geq 2$, there exists a $(r^2+1)$-dimensional family of simple (hence indecomposable)  Ulrich bundles of rank $r$ with Chern classes $c_1=rH$ and $c_2=\frac{dr^2+(2-d)r}{2}$. Moreover, they are $\mu$-semistable with respect to the polarization $-K_X=H=3e_0-\sum_{i=1}^{9-d}E_i$ and $\mu$-stable with respect to $H_{n}:=(n-3)e_0+H$ for $n\gg 0$.
\end{teo*}

Let us explain briefly how this paper is divided. In section $2$ we introduce the definitions and main properties that are going to be used throughout the paper; in particular, a brief account of ACM vector bundles, Ulrich bundles and Fano blow-up varieties is provided. In section $3$, we perform the construction of large families of simple ACM  vector bundles on all Fano blow ups of points in $\PP^n$. These families are constructed as the pullback of the kernel of certain surjective morphisms
$$\odi{\PP^n}(1)^b\longrightarrow\odi{\PP^n}(2)^a$$
with chosen properties. In the last section, we focus our attention on the $2$-dimensional case, where much more information is obtained. In the first subsection, we see that the ACM bundles we obtained by pullback are simple, Ulrich, and $\mu$-stable for a certain ample divisor $H_n$. Finally, in the last subsection, in the case of del Pezzo surfaces with very ample anticanonical divisors, we show that these families of vector bundles could also be obtained through Serre's correspondence from a general set of points on the surface.

\vskip 4mm \noindent {\bf Acknowledgements.} The authors would like to thank Costa and Mustopa for some fruitful discussions. We also  thank Mustopa for having shared with us their preprint \cite{CKM}.


\section{Preliminaries}
 Let us start fixing the notation that we are going to use through this paper. We are going to work with integral varieties over an algebraically closed field $k$ of zero characteristic. Given a smooth variety $X$ equipped with an ample line bundle $\odi{X}(1)$, the line bundle $\odi{X}(1)^{\otimes l}$ will be denoted by $\odi{X}(l)$ or $\odi{X}(l H)$. For any coherent sheaf $\shE$ on $X$ we are going to denote the twisted sheaf $\shE\otimes\odi{X}(l)$ by $\shE(l)$. As usual, $\Hl^i(X,\shE)$ (or simply $\Hl^i(\shE)$) stands for the cohomology groups and $\h^i(X,\shE)$ (or simply $\h^i(\shE)$) for their dimension. For a divisor $D$ on $X$, $\Hl^i(D)$ and $\h^i(D)$ abbreviate $\Hl^i(X,\odi{X}(D))$ and $\h^i(X,\odi{X}(D))$, respectively. When $\odi{X}(1)$ is very ample we will use the notation $\Hl^i_*(\shE)$ for the graded $\oplus_{l\in\mathbb{Z}} \Hl^0(X,\odi{X}(lH))$-module $\bigoplus_{l\in\mathbb{Z}} \Hl^i (X,\shE (l))$. $K_X$ will stand for the canonical class of $X$ and $\omega_X:=\odi{X}(K_X)$ for the canonical bundle.

\begin{defn}
We are going to say that a sheaf $\shE$ on $X$  is \emph{initialized} (with respect to $\odi{X}(1)$) if
$$
\Hl^0(X,\shE(-1))=0 \ \ \text{ but } \  \Hl^0(X,\shE)\neq 0.
$$
Notice that there always exists an integer $k$ such that $\shE_{init}:=\shE(k)$ is initialized.
\end{defn}

\begin{remark}\label{twist}
Given a rank $r$ vector bundle $\shF$ on a polarized variety $(X,\odi{X}(H))$, the first two Chern classes of $\shF$ are modified by twisting as follows:
\begin{itemize}
\item $c_1(\shF(lH))=c_1(\shF)+rlH$.
\item $c_2(\shF(lH))=c_2(\shF)+(r-1)lc_1(\shF)H+\binom{r}{2}l^2H^2$.
\end{itemize}
\end{remark}

\vskip 2mm
\noindent {\bf 2.1. ACM varieties and Ulrich bundles.}
In this subsection we are going to introduce the definitions and main properties of ACM varieties as well as those of ACM and Ulrich bundles needed later on.
\begin{defn}(cfr. \cite[Chapter I, Definition 1.2.2]{Mig}).
A closed subvariety $X\subseteq\PP^n_k$ is \emph{Arithmetically Cohen-Macaulay (ACM for short)} if its homogeneous coordinate ring $S_X$ is Cohen-Macaulay or, equivalently, $\di S_X=\depth \ S_X$.
\end{defn}

Notice that any zero-dimensional variety is ACM. For varieties of higher dimension we have the following characterization that will be used in this paper:
\begin{lemma}(cfr. \cite[ Chapter I, Lemma 1.2.3]{Mig}).\label{ACMcharacterization}
If $\di \ X\geq 1$, then $X\subseteq\PP^n_k$ is ACM if and only if $\Hl^i_*(\shI_X)=0$ for $1\leq i\leq \di X$.
\end{lemma}

\begin{defn}
Let $(X,\odi{X}(1))$ be a polarized variety. A coherent sheaf $\shE$ on $X$ is \emph{Arithmetically Cohen Macaulay} (ACM for short) if it is locally Cohen-Macaulay (i.e., $\depth \shE_x=\di \odi{X,x}$ for every point $x\in X$) and has no intermediate cohomology:
$$
\Hl^i_*(X,\shE)=0 \quad\quad \text{    for all $i=1, \ldots , \di X-1.$}
$$
\end{defn}
Notice that when $X$ is a smooth variety, which is going to be mainly our case, any coherent ACM sheaf on $X$ is locally free. For this reason we are going to speak uniquely of ACM bundles.

The notion of regularity with respect to a very ample line bundle is very classical. For our purposes, we will need to work in a slightly wider setting:
\begin{defn} (cfr. \cite[Definition 1.8.4]{Laz})
Let $X$ be a projective variety and $B$ an ample line bundle generated by its global sections. A coherent sheaf $\shF$ on $X$ is
\emph{m-regular with respect to B} if
$$
\Hl^i(X,\shF \otimes B^{(m-i)})=0 \text{ \ for $i>0$.}
$$

\end{defn}

\begin{teo}(cfr. \cite[Theorem 1.8.5]{Laz}) \label{lazarsfeld}
Let $\shF$ be an $m$-regular sheaf on $X$ with respect to $B$. Then for every $k\geq 0$:

\begin{enumerate}
\item[(i)]  $\shF\otimes B^{(m+k)}$ is generated by its global sections.
\item[(ii)] $\shF$ is $(m+k)$-regular with respect to $B$.

\end{enumerate}

\end{teo}

Let us now introduce the notion of a Ulrich sheaf.
\begin{defn}
Given a polarized $n$-dimensional variety $(X,\odi{X}(1))$, an ACM sheaf $\shE$ will be called an \emph{Ulrich sheaf} if $\h^0(\shE_{init})=\deg(X)\rk(\shE)$.
\end{defn}

\begin{remark}
In \cite{CKM}, initialized Ulrich bundles on hypersurfaces $X\subset \PP^n$ are called Clifford bundles. Moreover, if a hypersurface $X_{f}\subset\PP^n$ is defined by the equation $x_0^d-f(x_1,\dots,x_n)$ with $f$ an homogeneous polynomial $f\in k[x_1,...,x_n]$ of degree $d$ and $C_f$ denotes the Clifford algebra associated to f, which is defined to be the quotient of the free associated algebra $k\{\{y_0,...,y_n\}\}$ by the two-sided ideal generated by the relations $\{(x_1y_1+...+x_ny_n)^d=f(x_1,...,x_n) \mid x_1,...,x_n\in k \}$, they relate the problem of classifying all matrix representations of $C_f$ to the problem of studying the moduli space of Clifford bundles on $X_{f}$.

\end{remark}

When $\odi{X}(1)$ is very ample, we have the following result that justifies this definition:
\begin{teo}
Let $X\subseteq\PP^n$ be an integral subscheme and $\shE$  be an ACM sheaf on $X$. Then the minimal number of generators $m(\shE)$ of the $S_X$-module $\Hl^0_*(\shE)$ is bounded by
$$
m(\shE)\leq \deg (X)\rk(X).
$$
\end{teo}
\begin{proof}
See \cite[Theorem 3.1]{CH}.
\end{proof}

Therefore, since it is obvious that for an initialized sheaf $\shE$, $\h^0(\shE)\leq m(\shE)$, the minimal number of generators of Ulrich sheaves is
as large as possible. Modules attaining this upper bound were studied by Ulrich in \cite{Ulr}. A complete account is provided in \cite{ESW}. In particular, we have:
\begin{teo}(cfr. \cite[Proposition 2.1.]{ESW})\label{equivconditionsulrich}
Let $X\subseteq\PP^N$ be an $n$-dimensional ACM variety and $\shE$  be an initialized ACM coherent sheaf on $X$. The following conditions are equivalent:
\begin{enumerate}
\item[(i)] $\shE$ is Ulrich.
\item[(ii)] $\shE$ admits a linear $\odi{\PP^N}$-resolution of the form:

$$
0 \arr \odi{\PP^N}(-N+n)^{a_{N-n}}\arr \dots \arr \odi{\PP^N}(-1)^{a_1}\arr\odi{\PP^N}^{a_0}\arr\shE\arr 0.
$$
\noindent where $a_0=\de\shE$ and $a_i=\binom{N-n}{i}a_0$ for all $i$.
\item[(iii)] $\Hl^i(\shE(-i))=0$ for $i>0$ and $\Hl^i(\shE(-i-1))=0$ for $i<n$.
\item[(iv)] For some (resp. all) finite linear projections $\pi:X\arr\PP^n$, the sheaf $\pi_*\shE$ is the trivial sheaf $\odi{\PP^n}^t$ for some $t$.
\end{enumerate}
In particular, initialized Ulrich sheaves are $0$-regular and therefore they are globally generated.
\end{teo}

\noindent {\bf 2.2. Fano blow-ups of $\PP^n$.}
In this paper we are going to be interested in ACM bundles on blow-ups of points in $\PP^n$ with ample anticanonical divisor. We devote this subsection to recall their definition and main properties.

\begin{defn}(cfr. \cite[Chapter III, Definition 3.1]{Kol}).
A \emph{Fano variety} is defined to be a smooth $n$-variety $X$ whose anticanonical divisor $-K_X$ is
ample. Its degree is defined as $K_X^n$. If $-K_X$ is very ample, $X$ will be called a \emph{strong Fano variety}.
\end{defn}

\begin{remark}[Serre's duality for Fano varieties]
Let $X$ be a $n$-dimensional Fano variety with ample anticanonical divisor $H_X:=-K_X$. Given a vector bundle  $\shE$ on $X$, Serre's duality states:
$$
\Hl^i(X,\shE)\cong \Hl^{n-i}(X,\shE^{*}(-H_X))'.
$$
This remark will be used without further mention throughout the paper.
\end{remark}

In the following theorem we summarize the well-known results about the Picard group of the blow-up of $\PP^n$ at $s$ points and the intersection product of blow-ups needed in the sequel.
\begin{teo}(cfr. \cite[Chapter V, Proposition 4.8]{Harti}). \label{picblowups}
Let $\{p_1,\ldots ,p_s\}$  be a set of points in $\PP_k^n$ and let $\pi:X\rightarrow\PP_k^n$ be the blow-up of $\PP_k^n$ at
these points. Let $e_0\in Pic(X)$ be the pull-back of an hyperplane in $\PP^n$, let $e_i$ be the exceptional divisors (i.e., $\pi(e_i)=p_i$). Then:

\begin{enumerate}
\item[(i)] $Pic(X)\cong\mathbb{Z}^{s+1}$, generated by $e_0,e_1,\ldots ,e_s$.
\item[(ii)] The canonical class is $K_X=-(n+1)e_0+(n-1)\sum_{i=1}^s e_i$.
\item[(iii)] If $D\sim ae_0-\sum_{i=1}^s b_ie_i$, then $\chi (\OO_X(D))= {a+n\choose n}-\sum_{i=1}^s{b_i+n-1\choose n}$.
\item[(iv)] When $n=2$, the intersection pairing on the surface $X$ is given by $e_0^2=1$,$e_i^2=-1$, $e_0.e_i=0$ and $e_i.e_j=0$ for $i\neq j$.

\end{enumerate}

\end{teo}

In the particular two-dimensional case, Fano surfaces are called \emph{del Pezzo surfaces}. Their classification is classical and let us state it here for seek of completeness.

\begin{defn} (cfr. \cite{Dem})
A set of $s$ different points $\{p_1,\ldots ,p_s\}$ on $\PP_k^2$ with $s\leq 8$ is in \emph{general position} if no three of them lie on a line, no six of them lie on a conic and no eight of them lie on a cubic with a singularity at one of these points.
\end{defn}

\begin{teo}(cfr. \cite[Chapter IV, Theorems 24.3 and 24.4]{Man} and \cite[Prop. 8.1.9.]{Dol}).
Let $X$ be a del Pezzo surface of degree $d$. Then $1\leq d\leq 9$ and
\begin{enumerate}
\item[(i)] If $d=9$, then $X$ is isomorphic to $\PP_k^2$ (and $-K_{\PP^2_k}=3H_{\PP^2_k}$ gives the usual Veronese embedding in $\PP_k^9$).
\item[(ii)] If $d=8$, then $X$ is isomorphic to either $\PP_k^1\times\PP_k^1$ or to a blow-up of $\PP_k^2$ at one point.
\item[(iii)] If $7\geq d\geq 1$, then $X$ is isomorphic to a blow-up of $9-d$ closed points in general position.
\end{enumerate}
Conversely, any surface described under $(i),(ii),(iii)$ is a del Pezzo surface of the corresponding degree.
\end{teo}

\begin{lemma}(cfr. \cite[Prop. 3.4]{Kol}) \label{free}
Let $X$ be the blow-up of $\PP^2$ on $0\leq s\leq 8$ points in general position and let $K_X$ be the canonical divisor. Then:

\begin{enumerate}
\item[(i)] If $s\leq 6$, $-K_X$ is very ample and its global sections yield a closed embedding of $X$ in a projective space of dimension

$$ \di \Hl^0(X,\odi{X}(-K_X))-1=K_X^2=9-s.
$$

\item[(ii)] If $s=7$, $-K_X$ is ample and generated by its global sections.
\item[(iii)] if $s=8$, $-K_X$ is ample and $-2K_X$ is generated by its global sections.

\end{enumerate}
\end{lemma}

In the case of dimension $n\geq 3$ we are not allowed to blow up more than one point in $\PP^n$ in order to obtain a Fano variety. Indeed, we have:

\begin{teo}
Let $Z$ be a set of $s$ distinct points in $\PP^n$, $n\geq 3$, and let $X:=Bl_Z \PP^n\arr \PP^n$ be its blow-up. Then $X$ is Fano if and only if $s\leq 1$. Moreover, in this case $-K_X$ is very ample.
\end{teo}
\begin{proof}
The fact that the blow-up of $\PP^n$ at more than one point is not Fano is an immediate consequence of \cite[Theorem 1]{Bon}. On the other hand, it is obvious that $\odi{\PP^n}(n+1)$ is very ample. So let $X=Bl_{p}\PP^n$ be the blow-up of $\PP^n$ at one single point $p$. Its anticanonical divisor is $(n+1)e_0-(n-1)e_1$, which can be written as $(n-2)(e_0-e_1)$ plus $3e_0-e_1$. Since $e_0-e_1$ is clearly base-point free and $3e_0-e_1$ is very ample (which can be proven directly or appealing to the stronger result \cite[Theorem 1]{Cop}) we are done.
\end{proof}


\section{n-dimensional case}

This section contains the main result of this paper and its aim will be to exhibit (as far as we know) the first examples of $n$-dimensional varieties of wild representation type for arbitrary $n\ge 3$ (The case $n=2$ was already known. See, for instance, \cite{CH}, \cite{CKM} and \cite{PT}). More precisely, in this section we will construct large families of ACM vector bundles on Fano varieties of the form $X=Bl_Z \PP^n$ for $Z$ a finite set of $s$ distinct points on $\PP^n$.

We will keep the notation introduced in the previous section. In particular, $H$ will stand for the anticanonical ample divisor $-K_X$.   Recall from the previous section that in order to $X$ being Fano, we should assume that either $n=2$ and $Z$ is a set of up to $8$ points in general position or $n\geq 3$ and $s=0,1$.
It is well known that the embedding in $\PP^d$ via the anticanonical very ample line bundle $-K_X$ of a strong del Pezzo surface of degree $d$ is an ACM surface (cfr. \cite[Expos\'{e} V, Th\'{e}or\`{e}me 1]{Dem}); we are going to prove that the analogous result turns out to be true for Fano blow-ups of $\PP^n$, $n\ge 3$. In fact, we have

\begin{prop} Let $X=Bl_Z\PP^n$ be the blow-up of $\PP^n$, $n\ge 3$, at $s\le 1$ points and let us consider its embedding in $\PP^N_k$ through the very ample divisor $-K_X$. Then $X\subseteq\PP^N_k$ is an ACM variety.
\end{prop}
\begin{proof} Since it is well-known that Veronese embeddings are ACM, we can suppose that $X$ is the blow-up of $\PP^n$ at one single point. First of all, we are going to see that $\Hl^i_*(X,\odi{X})=0$ for $i=1,\dots, n-1$. To start with, notice that $\Hl^i(X,\odi{X})=\Hl^i(\PP^n,\pi_*(\odi{X}))=\Hl^i(\PP^n,\odi{\PP^n})=0$ for $i=1,\dots, n-1$. On the other hand, by Lemmas \ref{hi} and \ref{h1}, we have $\Hl^i(X,\odi{X}(tH))=0$ for $t>0$ and $i=1,\dots, n-1$. Finally, the vanishing of $\Hl^i(X,\odi{X}(tH))=0$ for $t<0$ is obtained by Serre's duality. So it would only remain to prove that $\Hl^1_*(\shI_{X|\PP^d})=0$, but this is immediate from the fact that $\Hl^1(X,\shI_{X|\PP^d})=0$ and that $\shI_{X|\PP^d}$ is $1$-regular.
\end{proof}

The ACM vector bundles $\shE$ on $X$ will be obtained as the kernel of certain surjective morphisms between $\odi{X}(e_0)^b$ and $\odi{X}(2e_0)^a$. Following notation from \cite{EHi}, let us consider $k$-vector spaces $A$ and $B$ of respective dimension $a$ and $b$. Set $V=\Hl^0(\PP^n,\odi{\PP^n}(1))$ and let $M=\Hom(B,A\otimes V)$ be the space of ($a\times b$)-matrices of linear forms. It is well-known that there exists a bijection between the elements of $m\in M$ and the morphisms $m:B\otimes\odi{\PP^n}\arr A\otimes\odi{\PP^n}(1)$. Taking the tensor with $\odi{\PP^n}(1)$ and considering global sections, we have morphisms $\Hl^0(m(1)):\Hl^0(\odi{\PP^n}(1)^b)\arr \Hl^0(\odi{\PP^n}(2)^a)$. The following result tells us under which conditions the aforementioned morphisms $m$ and $\Hl^0(m(1))$ are surjective:

\begin{prop}(\cite[Propoposition 4.1]{EHi})\label{hirschowitz}
For $a\geq 1$, $b\geq a+n$ and $2b\geq (n+2)a$, the set of elements from $m\in M$ such that $m$ and $\Hl^0(m(1))$ are surjective form a non-empty open set.
\end{prop}

Let us fix now the possible ranks $a$ and $b$, depending on the parity of $n\geq 2$, that we are going to deal with. If $n$ is even, we define:
\begin{equation}\label{even}
a=r \ \ \text{ and  } \ \ \ b=\frac{n+2}{2}r+c,
\end{equation}
\noindent for $r\geq 2$ and $c\in\{0,\dots,n/2-1\}$. If $n$ is odd,
\begin{equation}\label{odd}
a=2r \ \ \ \text{ and  } \ \ \ b=(n+2)r+c
\end{equation}
\noindent for $r\geq 1$ and $c\in\{0,\dots,n-1\}$. Notice that these values verify the conditions of Proposition \ref{hirschowitz}. So take an element
$m$ of the non-empty subset $U\subseteq M$ provided by Proposition \ref{hirschowitz} and consider the exact sequence
\begin{equation}\label{fam}
0\arr\shF\arr\odi{\PP^n}(1)^b\stackrel{m(1)}{\arr}\odi{\PP^n}(2)^a\arr 0.
\end{equation}

It is immediate to see that $\shF$, as a kernel of a surjective morphism of vector bundles, is a vector bundle of rank $rn/2+c$ in the even dimension case (resp. $nr+c$ in the odd dimension case). Let us consider now a finite set of $s$ distinct points $Z\subseteq\PP^n$ and the blow-up associated to these points $X:=Bl_Z(\PP^n)\stackrel{\pi }{\arr}\PP^n$. Pulling-back the exact sequence (\ref{fam}) we obtain he exact sequence:

\begin{equation}\label{fam2}
0\arr \pi^*\shF\arr\odi{X}(e_0)^b\stackrel{m(1)}{\arr}\odi{X}(2e_0)^a\arr 0.
\end{equation}

The first goal will be to show that $\shG:=\pi^*\shF$ is simple and therefore indecomposable. In order to do that, we are going to argue with the dual exact sequence

\begin{equation}\label{fam3}
0\arr\odi{X}(-2e_0)^a\stackrel{m(1)^{*}}{\arr}\odi{X}(-e_0)^b\arr \shG^* \arr 0.
\end{equation}

Notice that the morphism $f:=m(1)^*:\odi{X}(-2e_0)^a\longrightarrow \odi{X}(-e_0)^b$ appearing in the exact sequence  (\ref{fam3}) is a general element of the $k$-vector space $M:=\Hom(\odi{X}(-2e_0)^a,\odi{X}(-e_0)^{b})\cong k^{n+1}\otimes k^a\otimes k^{b}$ because $\Hom(\odi{X}(-2e_0),\odi{X}(-e_0))\cong\Hl^0(\odi{X}(e_0))\cong \Hl^0(\odi{\PP^n}(1))\cong k^{n+1}$. In other words, $f$ can be represented by a $b\times a$ matrix $A$ whose entries are elements of $\Hl^0(\odi{\PP^n}(1))$. Since $\Aut(\odi{X}(-e_0)^{b})\cong GL(b)$ and $\Aut(\odi{X}(-2e_0)^{a})\cong GL(a)$, the group $GL(b)\times GL(a)$ acts naturally on $M$ by

\begin{center}
\begin{tabular}{ c c l}
$GL(b)\times GL(a)\times M$ & $\longrightarrow$ & $M$ \\
$(g_1,g_2,A)$  & $\mapsto$ & $g_1^{-1}Ag_2$.\\
\end{tabular}
\end{center}

Moreover, for all $A\in M$ and $\lambda\in k^*$, $(\lambda Id_{b},\lambda Id_a)$ belongs to the stabilizer of $A$. Hence $\di_k Stab(A)\geq 1$. By \cite[Theorem 4]{Kac}, we have:

\begin{prop}\label{kac}
Let $M=k^{n+1}\otimes k^a\otimes k^{b}$ be endowed with the natural action of $GL(a)\times GL(b)$. If $a^2+b^2-(n+1)ab\leq 1$ then, for a general element $A\in M$,
$$
\di_k Stab(A)= 1.
$$
\end{prop}

The previous Proposition will allow us to conclude that the general vector bundle given by the exact sequence (\ref{fam3})
\begin{equation}\label{family4}
	  0 \arr \odi{X}(-2e_0)^a \stackrel{f}{\longrightarrow} \odi{X}(-e_0)^{b} \stackrel{g}{\longrightarrow} \shE:=\shG^* \arr 0
\end{equation}

\noindent is simple and, hence, indecomposable. More precisely, we have he following proposition which will be the key result for proving that all Fano blow ups of $\PP^n$ at a finite number of points are of wild representation type.

\begin{prop}\label{simple} Let $X=Bl_Z\PP^n$ be the blow-up of $\PP^n$ at a finite set of points, $n\geq 2$. If $n$ is even, we define:
$$
a=r \ \ \text{ and  } \ \ \ b=\frac{n+2}{2}r+c,
$$
\noindent for $r\geq 2$ and $c\in\{0,\dots,n/2-1\}$. If $n$ is odd, we set
$$
a=2r \ \ \ \text{ and  } \ \ \ b=(n+2)r+c,
$$
\noindent for $r\geq 1$ and $c\in\{0,\dots,n-1\}$. Let $\shF$ be the vector bundle obtained as the kernel of a general surjective morphism between $\odi{\PP^n}(1)^b$ and $\odi{\PP^n}(2)^a$:
\begin{equation}
0\arr\shF\arr\odi{\PP^n}(1)^b\stackrel{m(1)}{\arr}\odi{\PP^n}(2)^a\arr 0.
\end{equation}
Then, the vector bundles $\shE$ from pulling-back and dualizing $\shF$
\begin{equation}\label{mainsequence}
	  0 \arr \odi{X}(-2e_0)^a \stackrel{f}{\longrightarrow} \odi{X}(-e_0)^{b} \stackrel{g}{\longrightarrow} \shE:=(\pi^*\shF)^* \arr 0
\end{equation}
\noindent are simple.

\end{prop}
\begin{proof} First of all, notice that the values of $a$ and $b$ appearing in the statement of this Proposition verify the inequality of Proposition \ref{kac}.
We will now check that $\shE$ is simple.
Otherwise, there exists a non-trivial morphism $\phi:\shE\arr\shE$. Then we get a morphism $\overline{\phi}=\phi\circ g:\odi{X}(-e_0)^{b}\arr\shE$. Applying  $\Hom(\odi{X}(-e_0)^{b}, .)$ to the exact sequence(\ref{family4}) and taking into account
that $\Hom(\odi{X}(-e_0)^{b},\odi{X}(-2e_0)^a)=\Ext^1(\odi{X}(-e_0)^{b},\odi{X}(-2e_0)^a)=0$ we get $\Hom(\odi{X}(-e_0)^{b},\odi{X}(-e_0)^{b})\cong\Hom(\odi{X}(-e_0)^{b},\shE)$. Hence there is a non-trivial morphism $\widetilde{\phi}\in\Hom(\odi{X}(-e_0)^{b},\odi{X}(-e_0)^{b})$ induced by $\overline{\phi}$ and represented by a matrix $\mu Id \neq B\in GL(b)$ such that the following diagram commutes:
$$
\xymatrix{
0 \ar[r] &  \odi{X}(-2e_0)^a \ar[r]^f \ar[d]^{C} & \odi{X}(-e_0)^{b}  \ar[d]^{B} \ar[dr]^{\overline{\phi}} \ar[r]^{g}& \shE \ar[r] \ar[d]^{\phi} & 0\\
0 \ar[r] &  \odi{X}(-2e_0)^a \ar[r]^f & \odi{X}(-e_0)^{b} \ar[r]^g & \shE \ar[r] & 0.\\
}
$$
\noindent where $C\in GL(a)$ is the matrix associated to $\widetilde{\phi}_{|\odi{X}(-2e_0)^a}$. Then the pair $(\mu Id,\mu Id)\neq (C,B)$ belongs
to $Stab(A)$ and therefore $\di_k Stab(A)>1$ which produces the desired contradiction with Proposition \ref{kac}.

\end{proof}

Once the simplicity has been proved, the next goal will be to show that the bundles $\shE$ given by the exact sequence (\ref{family4}) are ACM. Since the proof in the case of dimension $2$ is slightly different and moreover one obtains a much stronger result, we postpone the discussion of this case until the next section and for the rest of the current one we only deal with $n\geq 3$. So, let $Z=\{p_1,\dots,p_s\}\subseteq\PP^n$ be a set of $s$ points on $\PP^n$ and let $X=Bl_Z\PP^n$ be the blow-up at this set. It is a well-known fact that for any integers $b_1, \cdots ,b_s\ge 0$ there is an isomorphism of cohomology groups $\Hl^i(X,\odi{X}(ae_0-\sum b_te_t))\cong\Hl^i(\PP^n,\shI_{W}(a))$ where $\shI_{W}$ is the ideal sheaf of the \emph{fat point} scheme $W=\sum_{t=1}^s b_tp_t$ (which is defined locally at $p_i$ by the ideal $\shI_{p_i}^{b_i}$ where $\shI_{p_i}$ is the maximal ideal of the local ring $\odi{\PP^n,p_i}$). This equivalence supplies us with the following standard result:
\begin{lemma}\label{hi}
Let $X=Bl_Z\PP^n$ be a blow-up at a set $Z$ of $s$ points in $\PP^n$. Assume that one of the following conditions holds:
\begin{itemize}\item[(i)] $2\leq i\leq n-1$, $a\in\mathbb{Z}$ and $b_t\geq 0$ for all $t$;
\item[(ii)] $i=0$ and $a<0$;  or \item[(iii)] $i=n$, $a\geq -n$.
\end{itemize} Then, we have $\Hl^i(X,\odi{X}(ae_0-\sum_{t=1}^sb_te_t))=0$.
\end{lemma}
\begin{proof}
(i) It is straightforward from the previous remark and the long exact sequence associated to
$$
0\arr\shI_W(a)\arr\odi{\PP^n}(a)\arr\odi{W}(a)\arr 0.
$$

(ii) If $a<0$, the divisor $ae_0-\sum_{t=1}^sb_te_t$ is not effective.

(iii) It follows from (ii) and Serre's duality.
\end{proof}

When we blow-up just one point we get easily the vanishing of some $\Hl^1$ groups needed later:
\begin{lemma}\label{h1}
Let $X=Bl_Z\PP^n$ be a blow-up at  $s=0,1$ points in $\PP^n$. If either $b=0$ or $a\geq b>0$ then $\Hl^1(X,\odi{X}(ae_0-be_1))=0$.
\end{lemma}
\begin{proof}
If $b=0$, then $\Hl^1(X,\odi{X}(ae_0))=\Hl^1(\PP^n,\odi{\PP^n}(a))=0$. On the other hand, if $a\geq b>0$, then a single point $p$ of multiplicity $b$ imposes ${b+n-1\choose n}$ independent conditions on hypersurfaces of degree $a$ and therefore from the exact sequence
$$
0\arr\shI_{bp}(a)\arr\odi{\PP^n}(a)\arr\odi{bp}(a)\arr 0
$$
\noindent follows immediately the vanishing of $\Hl^1(X,\odi{X}(ae_0-be_1))=0$.

\end{proof}

\begin{lemma}\label{surjective}
Let $X=Bl_Z\PP^n$ be a blow-up at a set $Z$ of $s\leq 1$ points in $\PP^n$ and let $$\odi{X}(e_0)^b\stackrel{f}{\arr}\odi{X}(2e_0)^a$$ be a morphism. Suppose that $f$ is surjective on global sections. Then for any $t\geq 1$, the induced morphism $\odi{X}(e_0+tH)^b\stackrel{f_t}{\arr}\odi{X}(2e_0+tH)^a$ is also surjective on global sections.
\end{lemma}
\begin{proof} If $s=0$ the result is obvious. Assume $s=1$, i.e. $Z=\{ p \}$.
Let us assume that the morphism $$\Hl^0(X,\odi{X}(e_0)^b)\stackrel{\Hl^0(f)}{\longrightarrow}\Hl^0(X,\odi{X}(2e_0)^a)$$ is surjective. Then, after taking the tensor product with $\Hl^0(\odi{X}(tH))$, the induced morphism
$$\begin{array}{c}
(\Hl^0(\odi{X}(e_0))\otimes\Hl^0(\odi{X}(tH)))^b\cong\Hl^0(\odi{X}(e_0)^b)\otimes\Hl^0(\odi{X}(tH))\stackrel{\Hl^0(f)}{\longrightarrow} \\
\Hl^0(\odi{X}(2e_0)^a)\otimes\Hl^0(\odi{X}(tH))\cong(\Hl^0(\odi{X}(2e_0))\otimes\Hl^0(\odi{X}(tH)))^a \end{array}
$$

\noindent is still surjective. For $t\geq 1$, let us consider the following commutative diagram:
$$
\xymatrix{
(\Hl^0(X,\odi{X}(e_0))\otimes\Hl^0(\odi{X}(tH)))^b \ar@{->>}[r] \ar[d] & (\Hl^0(X,\odi{X}(2e_0))\otimes\Hl^0(\odi{X}(tH)))^a  \ar[d] \\
 (\Hl^0(X,\odi{X}(e_0+tH)))^b\ar[r]^{\Hl^0(f_t)} &  (\Hl^0(X,\odi{X}(2e_0+tH)))^a  .\\
}
$$
So, in order to conclude the result, it will be enough to prove that the right vertical arrow on the previous diagram is surjective. From the above discussion, the previous diagram is equivalent to the following one:
$$
\xymatrix{
(\Hl^0(\odi{\PP^n}(1))\otimes\Hl^0(\PP^n,\shI_{W}(t(n+1))))^b \ar@{->>}[r] \ar[d] & (\Hl^0(\odi{\PP^n}(2))\otimes\Hl^0(\PP^n,\shI_{W}(t(n+1))))^a  \ar[d] \\
 (\Hl^0(\PP^n,\shI_{W}(t(n+1)+1)))^b\ar[r]^{\Hl^0(f_t)} &  (\Hl^0(\PP^n,\shI_{W}(t(n+1)+2)))^a  .\\
}
$$
\noindent where $\shI_W$ is the ideal sheaf of the fat point $W=t(n-1)p$. But, since by Lemmas \ref{hi} and \ref{h1}, we have  $\Hl^i(\shI_W(t(n+1)-2i))=0$ for $1\leq i\leq n$, $\shI_W(t(n+1))$ is $0$-regular with respect to $\odi{\PP^n}(2)$ and therefore the standard properties of regularity assure us that the right vertical arrow is surjective.

\end{proof}

Now we are ready to prove that the vector bundles that come out from the exact sequence (\ref{fam4}) are ACM in the case of dimension $n\geq3$. Remember that for these dimensions, we are allowed to blow-up a set $Z$ of $s=0,1$ points on $\PP^n$ in order to get a Fano variety $X=Bl_Z\PP^n$. Notice that the meaning of a blow-up of $s=0$ points is just a change of polarization of $\PP^n$, i.e. now we are considering $\PP^n$ with the very ample anticanonical divisor $(n+1)H$.

\begin{prop}\label{acm}
Let $X=Bl_Z\PP^n$ be a blow-up of $\PP^n$, $n\geq 3$, on $s=0,1$ points. Let $H:=-K_X$ be the ample anticanonical divisor. Then the vector bundle $\shE$ given by the exact sequence (\ref{mainsequence}) is ACM with respect to $H$.
\end{prop}

\begin{proof}
Let us consider the following pieces of the long exact sequence associated to the exact sequence (\ref{mainsequence}):
$$
\dots\arr \Hl^i(\odi{X}(tH-e_0)^b)\arr\Hl^i(\shE(tH))\arr\Hl^{i+1}(\odi{X}(tH-2e_0)^a)\arr\dots
$$
\noindent for $1\leq i\leq n-1$, and

$$
\dots\arr\Hl^n(\odi{X}(tH-e_0)^b)\arr\Hl^n(\shE(tH))\arr 0.
$$
Applying Lemma \ref{hi} we have that $\Hl^i(\shE(tH)))=0$ for all $i\geq 2$ and $t\geq 0$. On the other hand, for $t\geq 0$, by Lemma \ref{h1} we also get $\Hl^1(\shE(tH)=0$ for all $t\geq 0$. To see the remaining vanishings, let us consider the dual exact sequence (\ref{fam2})
\begin{equation}\label{fam4}
0\arr\shE^*\arr\odi{X}(e_0)^b\stackrel{m(1)}{\arr}\odi{X}(2e_0)^a\arr 0.
\end{equation}
Once again let us consider the following pieces of the long exact sequence
$$
\dots\arr \Hl^{i-1}(\odi{X}(tH+2e_0)^a)\arr\Hl^i(\shE^*(tH))\arr\Hl^i(\odi{X}(tH+e_0)^b)\arr\dots
$$
\noindent for $i\geq 2$ and $t\geq 0$. A new application of Lemmas \ref{hi} and \ref{h1} proves that $\Hl^{j}(\shE(sH))=\Hl^{n-j}(\shE^*((-s-1)H))=0$ for $0\leq j\leq n-2$ and $s\leq -1$. It only remains to show that $\Hl^{n-1}(\shE(tH))=\Hl^{1}(\shE^*((-t-1)H))=0$ for $t<0$. Notice that by Lemma \ref{h1} $\Hl^1(e_0+ tH)=0$ for all $t\geq 0$. Since we were in the case in which $\Hl^0(m(1))$ were surjective, by Lemma \ref{surjective}, we have that, in the exact sequence
$$
0\arr\Hl^0(\shE^*(tH))\arr\Hl^0(\odi{X}(e_0+tH)^b)\stackrel{f_t}{\arr}\Hl^0(\odi{X}(2e_0+tH)^a)\arr\Hl^1(\shE^*(tH))\arr 0,
$$

\noindent $f_t$ is surjective for all $t\geq 0$ and therefore we can conclude that $H^*\shE ^*(tH)=0$ for all $t\ge 0$ which proves what we want.

\end{proof}

We conclude the section gathering the previous results:
\begin{teo}\label{higheracm}
Let $X=Bl_Z \PP^n$ be a Fano blow-up of $\PP^n$, $n\geq 3$. Then, for any $r\geq n$, there exists a family of indecomposable ACM $r$-vector bundles of dimension $\frac{(n+2)n-4}{4}r^2-cr+c^2+1$ if $n$ is even or $((n+2)n-4)r^2-2cr+c^2+1$ if $n$ is odd. In particular, they are of wild representation type.
\end{teo}
\begin{proof} Let $X=Bl_Z \PP^n$ be a Fano blow-up of $\PP^n$, $n\geq 3$. For $r\geq 2$ if $n$ is even ($r\geq 1$ if $n$ is odd) let $a$ and $b$ be natural numbers defined as in (\ref{even}) and (\ref{odd}) and let $A$ and $B$ vector spaces of dimension respectively $a$ and $b$. We saw in Propositions \ref{simple} and \ref{acm} that a general element of the vector space $M=\Hom(B,A\otimes \Hl^0(\odi{\PP^n}(1)))$ provides with a simple ACM vector bundle $\shE$ on $X$ of rank $rn/2+c$ in the even dimension case (resp. $nr+c$ in the odd dimension case). Notice that in this way we sweep all the possible ranks bigger or equal than $n$. So it only remain to show that the dimension of the family is as it was stated. But this dimension can be computed as $\dim M - \dim\Aut(\odi{\PP^n}(1)^b)-\dim\Aut(\odi{\PP^n}(2)^a)+1$ which gives the announced result.
\end{proof}

\begin{remark} We want to point out that Theorem \ref{main} admits a natural generalization to a more general set up which will allow us to determine sufficient condition in order to assure that the $d$-uple embedding of $\PP^n$ is a variety of wild representation type. This will be discussed in a forthcoming paper (see \cite{MR1}) and in the Ph. D. Thesis of Pons-Llopis (see \cite{Pon}).
\end{remark}


\section{del Pezzo surfaces}
In this last section we focus our attention on $2$-dimensional Fano varieties $X$, i.e., on del Pezzo surfaces (excluding the quadric case $\PP^1\times\PP^1$) where much more information will be obtained. Recall that they are obtained as blow-ups of $\PP^2$ on $s\leq 8$ points in general position. From now on, when speaking of any del Pezzo surface we are excluding the quadric surface. $H$ will stand for the ample anticanonical divisor $-K_X$ on $X$. The aim is twofold: firstly, we are going to give an alternative construction of the family of ACM bundles on $X$ that will come with some extra information about the stability of the bundles. Secondly, we are going to show that these ACM bundles share a much stronger property: they are Ulrich bundles. In the final subsection, we are going to prove that, via the Serre's correspondence, the vector bundles could be obtain from finite general set of points on X verifying the Minimal Resolution Conjecture.

\vskip 4mm
\noindent {\bf 4.1. Construction of Ulrich bundles.}
The aim of this subsection is to recover, for any $r\geq 2$, the  $r^2+1$-dimensional family $\shE$ of rank $r$ ACM bundles on del Pezzo surfaces $X$ starting from rank $r$, $\mu $-stable vector bundles $\shH$ on $\PP^2$ with Chern classes $c_1(\shH) =0$ and $c_2(\shH) =r$.

One of the properties that we are going to discuss in this section is $\mu$-(semi)stability, let us recall the definition here. Given an $n$-dimensional variety $X$ and an ample divisor $H$ we say that $\shE$ is $\mu$-\emph{(semi)stable} with respect to $H$ if for every subsheaf $\shF$ of $\shE$ with $0<\rk\shF <\rk\shE$, $$\mu(\shF)<\mu(\shE) \quad (\text{resp. } \mu(\shF)\leq\mu(\shE)),$$ where the slope $\mu (\shE)$ of a vector bundle $\shE$ is defined as $$\mu(\shE):=\frac{c_1(\shE)H^{n-1}}{rk(\shE)}.$$ In \cite{Mar77} and \cite{Mar78}, M. Maruyama proved the existence of a projective coarse moduli scheme $M_{X,H}^{ss}(r;c_1,\dots,c_{min(r,n)})$ parameterizing certain equivalence classes of $\mu$-semistable vector bundles of rank $r$ and Chern classes $c_1,\dots,c_{min(r,n)}$. It contains an open (possibly empty) subscheme $M_{X,H}^{s}(r;c_1,\dots,c_{min(r,n)})$ parameterizing $\mu$-stable vector bundles.

Let us start proving  the existence of rank $r$, $\mu $-stable vector bundles $\shH$ on $\PP^2$ with Chern classes $c_1(\shH) =0$ and $c_2 (\shH) =r$.  In order to show that the moduli space $M^s_{\PP^2,\odi{\PP^2}(1)}(r;0,r)$ is non-empty we will use the following adapted result from \cite{DP}. Recall that the \emph{discriminant} of a vector bundle $\shE$ on $\PP^2$ of rank $r$ and Chern classes $c_1$, $c_2$ is defined as:
$$
\Delta(r,c_1,c_2)=\frac{1}{r}(c_2-\frac{(r-1)}{2r}c_1^2).
$$

\begin{teo}(cfr. \cite[Th\'{e}or\`{e}me B]{DP})
A sufficient and necessary condition for the existence of a $\mu$-stable vector bundle of rank $r$ and Chern classes $c_1\in r\mathbb{Z}$ and $c_2\in\mathbb{Z}$
on $\PP^2$ is that $\Delta(r,c_1,c_2)\geq 1$.
\end{teo}

\begin{cor}
The moduli space $M^s_{\PP^2,\odi{\PP^2}(1)}(r;0,r)$ is non-empty.
\end{cor}
\begin{proof}
Since in our case $\Delta=1$, we are done by the previous Theorem.
\end{proof}

Once we have checked the non-emptiness of $M^s_{\PP^2,\odi{\PP^2}(1)}(r;0,r)$ we can apply \cite[Proposition 4.3]{DMag} to assert that a generic element $\shH$ from $M_{\PP^2,\odi{\PP^2}(1)}(r;0,r)$ will have a resolution of the form
\begin{equation}\label{family}
	  0 \arr \bigoplus ^r \odi{\PP^2}(-2) \arr \bigoplus^{2r}\odi{\PP^2}(-1) \arr \shH \arr 0.
\end{equation}

Since this exact sequence is just the dual of the exact sequence constructed in (\ref{fam}), we see that we are recovering the family of vector bundles we were dealing with in section $3$.

The dimension of this family can be easily computed:

\begin{prop}
The dimension of the family of $\mu$-stable vector bundles from $M^s_{\PP^2,\odi{\PP^2}(1)}(r;0,r)$ with locally free resolution (\ref{family}) is $r^2+1$.
\end{prop}
\begin{proof}
By \cite[Theorem 4.4]{DMag}, $M^s_{\PP^2,\odi{\PP^2}(1)}(r;0,r)$ is an irreducible variety of dimension $r^2+1$. Since our family forms
an open subset of it, we are done.
\end{proof}

Let us take now a set of $0\leq s\leq 8$ points $Z=\{p_1,\dots,p_s\}$ in general position and let us consider the surface obtained by blowing up these points jointly with the canonical morphism to $\PP^2$, $$\pi:X=Bl_{Z}(\PP^2)\arr \PP^2.$$ We are following the notation introduced in Theorem \ref{picblowups}. Pulling back the vector bundles given by the exact sequence  (\ref{family}), we obtain the family of vector bundles from (\ref{fam4})

\begin{equation}\label{family2}
	  0 \arr \bigoplus ^r \odi{X}(-2e_0) \stackrel{f}{\longrightarrow} \bigoplus^{2r}\odi{X}(-e_0) \arr \shE:=\pi^*\shH \arr 0;
\end{equation}
\noindent If we twist it by the line bundle $\odi{X}(H)$ we obtain the following family of rank $r$
vector bundles:

\begin{equation}\label{family3}
	  0 \arr \bigoplus ^r \odi{X}(-2e_0+ H) \arr \bigoplus^{2r}\odi{X}(-e_0+H) \arr \shE(H) \arr 0;
\end{equation}

Their Chern classes can easily computed with the formulas given in Remark \ref{twist}: $$c_1(\shE(H))=rH \text{ and } c_2(\shE(H))=\frac{H^2r^2+(2-H^2)r}{2}.$$

We saw in Proposition \ref{simple} that these bundles were simple. Let us, however, provide an alternative proof that gives an stronger result. We will see that they are $\mu$-stable with respect a certain ample divisor and, therefore, simple. We are going to use the following result:

\begin{teo}(cfr. \cite[Theorem 1]{Nak})
Let $X,H$ be a polarized surface,  $\pi:X'\arr X$  the blow up of $X$ at $l$ distinct points $p_i$ and  $e_i$  the exceptional divisors. Let us define the divisor $H_n:=n\pi^*H-\sum_{i=1}^le_i$. Then for $n\gg 0$ there exists an open immersion $$\phi:M^s_{X,H}(r,c_1,c_2)\hookrightarrow M^s_{X',H_n}(r,\pi^{*} c_1,c_2)$$ defined by $\phi (\shF):= \pi^{*} (\shF)$ on closed points.
\end{teo}

\begin{cor}
The family of vector bundles on the blow up $\pi:X=Bl_{Z}(\PP^2)\arr \PP^2$ defined by the exact sequence
$$	  0 \arr \bigoplus ^r \odi{X}(-2e_0+ H) \arr \bigoplus^{2r}\odi{X}(-e_0+H) \arr \shE(H) \arr 0$$
 is $\mu$-stable with respect the ample divisor $ne_0-\sum e_i$ for $n$ big enough. In particular, they are simple, i.e., $\Hom(\shE(H),\shE(H))=k$.
\end{cor}

The last step will be to show, as in the higher dimensional case, that the bundles $\shE$ are ACM. In fact, much more will be provided in this case:
the twisted $\shE(H)$ bundles are initialized Ulrich bundles. For this, we need the following computations.

\begin{remark}[Riemann-Roch for vector bundles on a del Pezzo surface]
	Let $X$ be a del Pezzo surface. Since $X$ is a rational connected surface we have $\chi(\odi{X}) = 1$. In particular, Riemann-Roch formula for a vector bundle $\shE$ of rank $r$ has the form
	$$
	  \chi(\shE) = \frac{c_1(\shE)(c_1(\shE)-K_X)}{2} + r -c_2(\shE).
	$$
\end{remark}

\begin{remark}
The Euler characteristic of the involved vector bundles can be computed thanks to the Riemann-Roch formula:
\begin{equation}\label{xi1}
\chi(\odi{X}(-2e_0)(lH))=\frac{9-s}{2}l^2-\frac{3+s}{2}l,
\end{equation}

\begin{equation}\label{xi2}
\chi(\odi{X}(-e_0)(lH))=\frac{9-s}{2}l^2+\frac{3-s}{2}l,
\end{equation}

\begin{equation}\label{xi3}
\chi(\shE(lH))= 2r\chi(\odi{X}(-e_0)(lH)) -r\chi(\odi{X}(-2e_0)(lH))=\frac{9r-sr}{2}l^2+\frac{9r-sr}{2}l.
\end{equation}
\end{remark}

\begin{prop}\label{ulrich}
Let $X$ be a del Pezzo surface. The bundles $\shE(H)$ given by the exact sequence  (\ref{family3}) are initialized simple Ulrich bundles. Moreover, in the case of a blow-up of $\leq 7$ points, they are globally generated.
\end{prop}

\begin{proof}
First of all, notice that, since $\mu$-stability is preserved under duality, $\shE^*$ is a $\mu$-stable (with respect to $H_n$) vector bundle with $c_1(\shE^*)=0$, and therefore it does not have global sections: $\Hl^0(\shE^*)=\Hl^2(\shE(-H))=0$. In particular, $\Hl^2(\shE(tH))=0$, for all $t\geq -1$.
On the other hand, since $\Hl^2(\odi{X}(-2e_0))=\Hl^0(\odi{X}(2e_0-H))=0$ and $\h^1(\odi{X}(-e_0))=-\chi(\odi{X}(-e_0))=0$ we obtain from the long
exact sequence of cohomology associated to (\ref{family2}) that $\Hl^1(\shE)=0$. Since $\chi(\shE)=0$, we also conclude that $\Hl^0(\shE)=0$ and therefore $\Hl^0(\shE(tH))=0$ for all $t\leq 0$. Moreover, since we also have that $\chi(\shE(-H))=0$, we obtain that $\Hl^1(\shE(-H))=0$.

Now, it is well-known that $\h^0(\odi{X}(H))=H^2+1>0$ (see for instance \cite[Corollary 3.2.5]{Kol}) and therefore there exists an exact sequence:

$$
 0 \arr \odi{X}(-H) \arr \odi{X} \arr \odi{H} \arr 0.
$$

If we tensor it by $\shE$ and we consider the long exact sequence associated to it we see that

$$
0=\Hl^0(\shE) \arr \Hl^0(\shE_{\mid H}) \arr \Hl^1(\shE(-H))=0.
$$

\noindent This shows that $\Hl^0(\shE_{\mid H}(-tH))=0$ for all $t\geq 0$. Then we can use this last fact jointly with the long exact sequence associated to
$$
 0 \arr \shE(-(t+1)H) \arr \shE(-tH) \arr \shE_{\mid H}(-tH) \arr 0
$$
to show inductively that $\Hl^1(\shE(-tH))=0$ for all $t\geq 0$.

In order to complete the proof we need to consider now two different cases:
\begin{itemize}
\item \emph{$X$ is the blow-up of $s\leq 7$ points on $\PP^2$ in general position}. In this case, by Lemma \ref{free}, $H$ is ample and generated by its global sections. Since we have just seen that $\shE(H)$ is $0$-regular with respect to $H$ we can conclude by Theorem \ref{lazarsfeld} that $\shE(H)$ is ACM and globally generated. Moreover, $\h^0(\shE(H))=\chi(\shE(H))=(9-s)r=H^2r$, i.e., $\shE(H)$ is an Ulrich bundle.
\item \emph{$X$ is the blow-up of $8$ points on $\PP^2$ in general position}. In this case, the argument is slightly more involved, since we can use Theorem \ref{lazarsfeld} only with respect to $2H$, which is ample and globally generated. First of all, since the points are in general position, $\Hl^0(\odi{X}(-e_0+H))=0$ and from the exact sequence (\ref{family3}) we get the following exact sequence:
    $$
    0 \arr \Hl^0(\shE(H)) \arr \oplus^r \Hl^1(\odi{X}(-2e_0+H)) \arr \oplus^{2r}\Hl^1(\odi{X}(-e_0+H)) \arr \Hl^1(\shE(H)) \arr 0.
    $$
    From this sequence and the fact that $\h^1(\odi{X}(-2e_0+H))=-\chi(\odi{X}(-2e_0+H))=5$ and  $\h^1(\odi{X}(-e_0+H))=-\chi(\odi{X}(-e_0+H))=2$ we are forced to conclude that $\h^0(\shE(H))=r$ and $\Hl^1(\shE(H))=0$. Now, from what we have gathered up to now, we can affirm that $\shE(H)$ is $1$-regular with respect to $2H$ and therefore, by Theorem \ref{lazarsfeld}, $\Hl^1(\shE(H+2tH))=0$ for all $t\geq 0$. In order to deal with the cancelation of the remaining groups of cohomology, it will be enough to show that $\shE(2H)$ is $1$-regular with respect to $2H$, i.e., it remains to show that $\Hl^1(\shE(2H))=0$. In order to do this consider the exact sequence (the cancelation of $\Hl^0(\odi{X}(-e_0+2H))$ is due to the fact that the points are in general position):
    $$
    0 \arr \Hl^0(\shE(2H)) \arr \oplus^r \Hl^1(\odi{X}(-2e_0+2H)) \arr \oplus^{2r}\Hl^1(\odi{X}(-e_0+2H)) \arr \Hl^1(\shE(2H)) \arr 0.
    $$
    \noindent Once again, we control the dimension of these vector spaces:  $\h^1(\oplus^r\odi{X}(-2e_0+2H))=-r\chi(\odi{X}(-2e_0+2H))=9r$ and
    $\h^1(\oplus^{2r}\odi{X}(-e_0+2H))=-2r\chi(\odi{X}(-e_0+2H))=6r$ Therefore we are forced to have $\h^0(\shE(2H))=3r$ and $\Hl^1(\shE(2H))=0$. Notice that in this case $\shE(3H)$ is globally generated.
\end{itemize}

\end{proof}

Given a del Pezzo surface $X$, we have just seen that the bundles given by the exact sequence (\ref{family3}) were $\mu$-stable with respect to the ample divisor $H_n:= ne_0 -\sum e_i$ for $n$ big enough. Unfortunately, the proof did not provide an effective value of $n$. However we are going to prove at least that they are $\mu$-semistable with respect to the anticanonical divisor $H=H_3=3e_0-\sum e_i$. The main tool will be the classification of vector bundles on elliptic curves performed in  \cite{Ati}:

\begin{prop}
Let $X$ be a del Pezzo surface of degree $d$. Then a general bundle $\shE(H)$ given by the exact sequence  (\ref{family3}) is $\mu$-semistable.
\end{prop}

\begin{proof}
We follow the structure of the proof given by the  case of the cubic surface in \cite[Proposition 5.2]{CH}. We saw in Proposition \ref{ulrich} that these vector bundles $\shE(H)$ were initialized and Ulrich. Moreover, we know that we can take a smooth elliptic curve $H$ as a representative of the anticanonical divisor class (see for instance \cite[III,Theorem 1]{Dem}). From the exact sequence
$$
 0 \arr \shE(-H) \arr \shE \arr \shE_{|H} \arr 0
$$
is deduced that $\shE_{|H}(H)$ is also initialized of degree $dr$ and $\h^0(\shE_{|H}(H))=dr$. By \cite[Theorem 7]{Ati}, $\shE_{|H}(H)=\oplus\shE_{r_i,d_i}$ with $\shE_{r_i,d_i}$ a vector bundle of rank $r_i$ and degree $d_i$. Since $\h^0(\shE_{r_i,d_i}(-H))=0$, Atiyah's classification forces that $d_i\leq dr_i$. It follows that we have equality and $\shE_{|H}(H)$ decomposes as direct sum of $\mu$-semistable bundles of the same slope. Thus, $\shE_{|H}(H)$ is also  $\mu$-semistable of slope $d$. Therefore it is straightforward to conclude that $\shE(H)$ is also $\mu$-semistable.
\end{proof}

Summing up, we get the following result:
\begin{teo}\label{main}
Let $X$ be a del Pezzo surface  of degree $d$. Then for any $r\geq 2$ there exists a family of dimension $r^2+1$ of simple initialized Ulrich bundles of rank $r$ with Chern classes $c_1=rH$ and $c_2=\frac{dr^2+(2-d)r}{2}$. Moreover, they are $\mu$-semistable with respect the polarization $H=3e_0-\sum_{i=1}^{9-d}E_i$ and $\mu$-stable with respect $H_{n}:=(n-3)e_0+H$ for $n\gg 0$. In particular, del Pezzo surfaces are of wild representation type.
\end{teo}

\begin{remark}
As in the case of surfaces on $\PP^3$, studied in \cite[Lemma 2.16 and Proposition 3.5]{CKM}, it is possible to see that if $\shE$ is a rank $r$  initialized Ulrich bundle on a strong del Pezzo surface $X$ of degree $d$, then $$ \begin{array}{l} c_1(\shE)H=dr, \\ c_2(\shE)=\frac{c_{1}^2+(2-d)r}{2}, \text{ and} \\
\frac{(d-2)r^2+(2-d)r)}{2} \leq c_2 \leq\frac{dr^2+(2-d)r}{2}.
\end{array}$$
Therefore, the Ulrich bundles constructed in Theorem \ref{main} are extremal with respect to the second Chern class.
\end{remark}

\begin{remark}
\begin{itemize}
\item The existence of rank $2$ Ulrich bundles on arbitrary del Pezzo surfaces were established by Eisenbud, Schreyer and Weyman in \cite[Corollary 6.5]{ESW}.
\item On smooth cubic surfaces, the existence of $(r^2+1)$-dimensional families of rank $r$ Ulrich bundles with the aforementioned Chern classes was announced by Casanellas and Hartshorne in \cite{CH} (see also \cite{CKM}).
\end{itemize}
\end{remark}


\vskip 4mm
\noindent {\bf 4.2.Serre's correspondence.}
In this last subsection we are going to pay attention to the case of strong del Pezzo surfaces $X$. In this case, the very ample divisor $-K_X$ provides an embedding $X\subseteq\PP^d$, with $d=K_X^2$. Let $R:=k[X_0,\dots, X_d]$ be the graded polynomial ring associated to $\PP^d$. We are going to show that the $(r^2+1)$-dimensional family of rank $r$ Ulrich bundles given in Theorem \ref{main} could also be obtained through a version of Serre's correspondence from a general set of $c_2(\shE(H))=\frac{dr^2+ (2-d)r}{2}$ points on $X$.

\vskip 2mm
We are going to use the fact that  $\frac{dr^2+ (2-d)r}{2}$ general points on a del Pezzo surface $X\subseteq\PP^d$ satisfy the Minimal Resolution Conjecture stated by Musta\c{t}\u{a} in \cite{Mus}. More precisely, we have the following result:

\begin{teo}( cfr.\cite[Theorem 3.7]{MRP})\label{gorenstein}
Let $X\subseteq\PP^d$ be a strong del Pezzo surface of degree $d$ embedded in $\PP^d$ by its very ample anticanonical divisor. Let $Z_{m(r)}\subset X$  be a general set of $m(r)=1/2(dr^2+ (2-d)r)$ points, $r\geq 2$. Then the minimal graded free resolution (as a $R$-module) of the saturated ideal of $Z_{m(r)}$ in $X$ has the following form:
\begin{equation}\label{res}
0\longrightarrow R(-r-d)^{r-1}\longrightarrow
R(-r-d+2)^{\gamma_{d-1,r-1}}\longrightarrow\ldots
\end{equation}
$$
\longrightarrow R(-r-1)^{\gamma_{2,r-1}}\longrightarrow
R(-r)^{(d-1)r+1}\longrightarrow I_{Z_{m(r)}|X}\longrightarrow 0
$$
with
$$\gamma_{i,r-1}
=\sum_{l=0}^{1}(-1)^{l}\binom{d-l-1}{i-l}\Delta^{l+1}P_X(r+l)-\binom{d}{i}(m(r)-P_X(r-1)).
$$

\end{teo}

\begin{teo}(Serre's correspondence)\label{serre}
Let $X\subseteq\PP^d$ be a strong del Pezzo surface of degree $d$.
\begin{enumerate}
\item[(i)] If $\shE(H)$ is an Ulrich bundle of rank $r\geq 2$ given by the exact sequence (\ref{family3}), then there is an
exact sequence
$$
0 \arr \odi{X}^{r-1} \arr \shE(H) \arr \shI_{Z|X}(rH) \arr 0
$$

\noindent where $Z$ is a zero-dimensional scheme of degree $m(r)=c_2(\shE(H))=1/2(dr^2+ (2-d)r)$ and $\h^0(\shI_{Z|X}(r-1))=0$.

\item[(ii)] Conversely, for general sets $Z$ of $m(r)=1/2(dr^2+ (2-d)r)$ points on $X$, $r \geq 2$, we recover the initialized Ulrich bundles given by the exact sequence (\ref{family3}) as an extension of $\shI_{Z,X}(rH)$ by $\odi{X}^{r-1}$.
\end{enumerate}

\end{teo}
\begin{proof}
(i) As $\shE(H)$ is globally generated, $r-1$ general global sections define an exact sequence of the form

$$
0 \arr \odi{X}^{r-1} \arr \shE(H) \arr \shI_{Z|X}(D) \arr 0
$$
\noindent where $D=c_1(\shE(H))=rH$ is a divisor on $X$ and $Z$ is a zero-dimensional scheme of length $$c_2(\shE(H))=\frac{dr^2+ (2-d)r}{2}.$$
Moreover, since $\shE(H)$ is initialized, $\h^0(\shI_{Z|X}(r-1))=0$.

(ii) This is a classical argument so we are going to be brief. We follow the lines of the argument given for  the cubic case in \cite[Theorem 4.4]{CH}. Let $Z$ be a general set of points of cardinality $m(r)$ with the minimal free resolution of (\ref{res}). Let us denote by $S_X$ and $S_Z$ the homogeneous coordinate ring of $X$ and $Z$. It is well-known that for Arithmetically Cohen-Macaulay varieties, there exists a bijection between ACM bundles on $X$ and Maximal Cohen Macaulay (MCM from now on) graded $S_X$-modules sending $\shE$ to $\Hl^0_*(\shE)$. From the exact sequence
$$
0 \arr I_{Z|X} \arr S_X \arr S_Z \arr 0
$$
\noindent we get $\Ext^1(I_{Z|X},S_X(-1))\cong\Ext^2(S_Z,S_X(-1))\cong K_Z$ where $K_Z$ denotes the canonical module of $S_Z$ (the last isomorphism is due to the fact that $S_X(-1)$ is the canonical module of $X$ and the codimension of $Z$ in $X$ is $2$). Dualizing the exact sequence (\ref{res}), we obtain a minimal resolution of $K_Z$:
$$
\dots \arr R(r-3)^{\gamma_{d-1,r-1}} \arr R(r-1)^{r-1} \arr K_Z \arr 0.
$$
\noindent This shows that $K_Z$ is generated in degree $1-r$ by $r-1$ elements. These generators provide an extension
\begin{equation}\label{aa}
0\arr S_X^{r-1} \arr F \arr I_{Z|X}(r) \arr 0
\end{equation}
via the isomorphism $K_Z\cong\Ext^1(I_{Z|X},S_X(-1))$. $F$ turns out to be a MCM module because $\Ext^1(F,K_X)=0$
(this last cancelation follows by applying $\Hom_{S_X}(-,K_X)$ to (\ref{aa})).
If we sheafify the exact sequence (\ref{aa}) we obtain the sequence

$$
0\arr \odi{X}^{r-1} \arr \widetilde{F} \arr I_{Z|X}(r) \arr 0
$$

\noindent where $\widetilde{F}$ is an ACM vector bundle on $X$. From (\ref{res}) it can be seen that $\Hl^0(\shI_{Z|X}(r-1))=0$ and $\h^0(\shI_{Z|X}(r))=(d-1)r+1$
and therefore $\widetilde{F}$ will be a Ulrich bundle (i.e., $\h^0(\widetilde{F})=dr$) and initialized. By Theorem \ref{equivconditionsulrich}, $\widetilde{F}$ will be globally generated.

It only remains to show that for a generic choice of $Z_{m(r)}\subset X$, the associated bundle $\shF:=\widetilde{F}$ just constructed belongs to the family (\ref{family3}). Since $\shF$ is an initialized Ulrich bundle of rank $r$ with the expected Chern classes, the problem boils down to a dimension counting. We need to show that the dimension of the family of bundles obtained through this construction from a general set $Z_{m(r)}$ is $r^2+1$.
Since this dimension is given by the formula $\di Hilb^{m(r)}(X)- \di Grass(\h^0(\shF),r-1)$, an easy computation taking into account that $\di Hilb^{m(r)}(X)=2m(r)$ and $\di Grass(\h^0(\shF),r-1)=(r-1)(dr-r+1)$ gives the desired result.

\end{proof}

\end{document}